\documentclass{article}
\usepackage{amsmath}
\usepackage[english]{babel}
\usepackage[utf8]{inputenc}
\usepackage{graphicx}
\usepackage{soul}
\usepackage{epstopdf}
\usepackage[bottom=2.5cm,left=3cm,right=3cm,top=3cm]{geometry}

\begin{document}


\begin{center}
\section*{\LARGE Spectrum for Discrete Closed Chain of Contours with Cluster Movement}

\textit{Myshkis P.A.$^{1}$, Tatashev A.G.$^{1}$, Yashina M.V.$^{1}$}\\
\textit{$^{1}~$Moscow Automobile and Road Construction State Technical University (MADI)}\\
e-mail:yash-marina@yandex.ru
\end{center}

\large 
{\bf Abstract.} This paper studies a discrete dynamical system  belonging to the class of the networks introduced by A.P.~Buslaev. The systems contains a finite set of contours. In any contour, there are cells and a group of particles. This group is called a cluster. The number of these particles is even. The particles are in adjacent cells and move simultaneously. For each contour, there are two adjacent contours. These contours are the contour on the left and the contour on the right. There is a common point for two adjacent contours. These common points are called nodes. At each discrete moment, the particles of a cluster move onto a cell forward. Delays in the cluster movement occur due to that particles of two clusters may not cross the same node simultaneously. The main problem is to study limit cycles of the considered dynamical systems and the set of realized velocities of clusters taking into account the delays.

\section{Introduction}
A.P. Buslaev has developed a class of dynamical systems, which were later called the contour networks, with aim to elaborate network traffic mathematical models such that analytical results may be obtained [1].

A contour network contains a system of closed sequences of cells called contours. There are common points of adjacent cells. These points are called nodes. Discrete and continuous versions of contour networks are considered. In the case of a discrete version, each  contour contains cells. At any discrete moment, each cell is vacant or occupied by a particle. Two types of particle movement are considered.
These types are  the individual movement and the cluster movement.

For individual movement, if the cell ahead of a particle is vacant, then the particle moves onto a cell forward, and the cell is in the cell ahead at next step. This is consistent with the rule of movement for a simple version of the Nagel~-- Schreckenberg traffic model [2]. 
For the traffic models with these rules of movement, analytical results were obtained under the assumption that the particles move along a one-dimensional closed or infinite lattice, for example, in [6--8]. A traffic model with a similar type of movement but for two-dimensional toroidal lattice was studied in [9--12]. 

The concept of cluster movement in mathematical models of traffic was introduced in  [13]. In the case of this type of movement, particles form gropes such that the particles of the same group occupy adjacent cells and move simultaneously. 

For continuous contour networks, a cluster is a segment moving along a contour with constant velocity if there are no delays. Delays of clusters (particles) occur at nodes and are due to that more than one cluster cannot move through the same node simultaneously as in the case of discrete contour networks. 

 Contour networks were mainly considered such that transitions of particles (clusters) from a contour to another contour are not allowed. 

Analytical results are obtained for two-contours systems with one [14, 15] or two [16, 17] nodes, contour networks with one-dimensional periodic contour networks (closed [1, 18--22] and open [23, 24] chains of contours), contour networks with two-mensional toroidal [1, 25, 26]
or rectangular structure [26].

This paper studies a discrete closed chain of contours with the same number cells in each contour and the same number of particles in each contour. For any contour, there is a node common with each of two adjacent contours. The nodes divide the contours into two parts of the same length. The nodes are located between two cells at each contour. A delay in the movement of a cluster if the cluster comes to a node when the cluster of the adjacent contour passes through this node. If two clusters comes to the same node simultaneously, a competition occurs. The competition is resolved according to the left-priority rule, i.e., the cluster located to the left of the node passes the node first. The spectrum of the system is studied. The concept of the spectrum, i.e., the set of limit cycles of the system with related average velocities of clusters taking into account delays was introduced in [20]. The paper [20] studies almost the same system as the present paper. The paper [20] studies the spectrum of the system with the number of contours equal to 2 or 3. In [19], analogous results were obtained for the continuous version of the closed chain with two or three contours. The papers [18, 21] study the spectrum of a discrete closed chain with different competition resolution rule  such that, in each contour of the chain, there are two cells and a particle (binary chain of contours). In [21], it is proved that, for the closed chain of contours with prescribed number cells in the contour, the number of contours can be prescribed such that there exist limit cycles with delays of particles, i.e., for any small density of particles,  limit cycles with delays may exist. 
  In [22], an analogous was proved for a continuous closed chain of contours. 
  
 This paper studies the spectrum of discrete closed chain with an arbitrary number of contours. There is a cluster at each contour. An algorithm is obtained to find limit cycles. This paper proves the hypothesis formulated in [20] for the spectrum of closed chain with an arbitrary number of contours.

Section 2 describes the studied system. In Section~3, it is proved that a necessary and sufficient condition for the system to result in the state of collapse (from a moment, all clusters do not move at present and in the future) is that the length cluster be not greater than a half of the length of contour.  In Section 4, conditions have been found for that there exist limit cycles such that the system is not in the state of collapse but delays occur. In Section~5, a sufficient condition for the system self-organization, i.e., for that, from any initial state, the system results in the state such that all clusters move without delays at present and in the future. Section~6 proves a theorem regarding the set of realized velocities of clusters.  

\section{Description of system}
Let the system contain $N$ {\it contours} with indexes $i=0,1,\dots,N-1.$ There are $2m$ {\it cells}, $m\ge 1,$ and a {\it cluster} containing $l<2m$ {\it particles.} The particles are adjacent cells and move simultaneously. At any discrete moment $t=0,1,2,\dots$ a cluster moves in the prescribed direction. The cells are numbered in the direction of movement. The indexes of the cells are $0,1,\dots,2m-1.$ For the contour $i,$ the contour adjacent to the left is the contour $i-1$ (subtraction modulo $N),$ and the contour $i+1$ adjacent to the right is the contour $i+1$ (subtraction modulo $N),$ $i=0,1,\dots,N-1.$ If, at present moment, the leading particle of the cluster of the contour $i$ (the {\it cluster $i)$} is in the cell $j,$ then we say that the contour $i$ is in the state $j,$ $i=0,1,\dots,N-1,$ $j=0,1,\dots,2m-1.$ In this case, the other particles of the cluster $i$ are in the cells $j-1,\dots,j-m+1$ (addition and subtraction modulo $m).$ For the contours $i,$ $i+1$ (addition modulo $N),$  there is a common node called the node $(i,i+1)$ and is between the cells 0, 1 of the contour $i$ and between the cells $m,$ $m+1$ of the contour $i+1,$ $i=0,1,\dots,N-1.$ We say that the cluster $i$ {\it occupies the node $(i,i+1),$} 
$i=0,1,\dots,N-1,$ if
the contour $i$ is in one of the states $1,\dots,l-1.$ 
We say that the cluster $i$ {\it occupies the node $(i-1,i),$} 
$i=0,1,\dots,N-1,$ if the contour $i$ is in one of the states $m+1,\dots,m+l-1$ (addition modulo $N).$
Suppose, at time $t,$ the contour is in the state $j.$ Then, at time $t+1,$ the contour is in the state $j+1$ (addition modulo $2m),$ $j=0,1,\dots,2m-1,$ if no delay occurs. We say that the cluster $i,$ {\it is at the node $(i,i+1)$} if the contour $i,$ $i=0,1,\dots,N-1,$ in the state 0. We say that the cluster $i,$ {\it is at the node $(i,i-1)$} if the contour $i,$ $i=0,1,\dots,N-1,$ is in the state $m.$ A delay in the movement of a cluster occurs if, at time $t,$ this cluster is at a node, and the cluster of the adjacent contour occupies this node. In this case, the cluster being at time $t$ at a node, is, at time $t+1,$ will be at this node. If the clusters  $i$ and $i+1$ (addition modulo $N)$ are at the same node simultaneously, then the cluster $i$ passes the node first (the {\it left-priority competition resolution rule}), $i=0,1,\dots,N-1.$  

The state of the system at time $t$ is the vector
$$(x_0(t),x_1(t),\dots,x_N(t)),$$  
where $x_i(t)$ is the state of the contour $i$ at time $t,$ $i=0,1,\dots,N-1,$ $t=0,1,2,\dots$ 

A state is{\it admissible} if, for this state, no node is occupied by two clusters.  

Evidently, at any time, no node is occupied more than one cluster. 

Suppose
$$\Delta_i(t)=x_{i+1}(t)-x_i(t),\eqno(1)$$
(subtraction modulo $2m,$ and sum in the index modulo $N).$

\section{Average velocity of particles. State of free movement. State of collapse}

The considered system is deterministic, and the number of the system states is finite. Since that a state repeated, a sequence of states is repeated periodically. Denote by $T$ the period. The value $T$ is the duration of limit cycle of dynamical system. 

Let there are $S_i(t)$ transitions of the cluster $i$ in the time interval $(0,t).$ The limit
$$v_i=\lim\limits_{t\to \infty}\frac{S_i(t)}{t}$$    
is called {\it the average velocity of the cluster} $i,$ $i=0,1,\dots,N-1.$ Obviously, that this limit exists and is equal to the ratio of the number of transitions of the cluster $i$ during the limit cycle to the period $T.$  

We shall prove that the average velocity of any cluster is the same. Denote the average velocity of clusters by $v.$    

We say that, at time $t_0,$ the system {\it is in the state of free movement} if at any time $t\ge t_0,$ any cluster moves without delays. If the system results in the state of free movement, then $v=1.$

We say that, at time $t_0,$ the system {\it is in the state of collapse} if, at time $t\ge t_0,$ no cluster moves. If the system results in the state of free collapse, then $v=0.$ 

\section{Necessary and sufficient condition for that the system results in the state of collapse}

In this section, it is proved that a necessary and sufficient condition for the system to result in the state of collapse (from a moment, all clusters do not move at present and in the future) is that the length cluster be not greater than a half of the length of contour.  
\vskip 3pt
{\bf Theorem 1.} {\it If the condition
$$l\le m\eqno(2)$$
 holds, then there exist no states of collapse.  

If the condition (2) does not hold that the system results in the state of collapse from any initial state.  
\vskip 3pt
Proof.} If (2) holds, then any cluster occupies at least one node. There are  $N$ clusters and $N$ nodes. No node can be occupied more than one cluster simultaneously. Hence no cluster an occupy two nodes simultaneously. Therefore, if we assign a cluster to the node occupied by this cluster, then we get a one-to-one correspondence between the set of clusters and the set of clusters and the set of nodes. However, if a leading particle of a cluster passes a node, then, at next step, either the cluster will occupy two nodes (if $l\ge m+2)$ or a node, which is occupied by cluster will be vacant if $l=m+1.$ Thus, if (2) holds, then no cluster cross a node.

\section{Limit cycles with delays}
Assume that $l\le 1/2,$ and hence according to Theory~1 the system results in the state of collapse. 

A delay of the cluster $i+1$ (addition modulo $N)$ at the node $(i,i+1)$  is called a {\it delay of the first type,} and a delay of the cluster $i$ at the node $(i,i+1)$ is called a  {\it delay of the second type,} $i=0,1,\dots,N-1.$ The duration of the delay of the first type delay is not greater than $l,$ and the duration of the second type is not greater than $l-1.$  
\vskip 3pt
{\bf Lemma 1.} {\it If, at time $t,$ a delay of the first type for the cluster $i+1$ ends, then 
$\Delta_i(t)=m-l,$ $i=0,1,\dots,N-1.$   
\vskip 3pt
Proof.} If, at time $t,$ a delay of the first type for the cluster $i+1$ ends, then   
$x_i(t)=l,$ $x_{i+1}(t)=m.$ From this and (1), we get the lemma. 
\vskip 3pt
{\bf Lemma 2.} {\it If, at time $t,$ a delay of the second type for the cluster $i,$  then $\Delta_i(t)=m+l,$ $i=0,1,\dots,N-1.$}
\vskip 3pt
The proof of Lemma 2 is the same as the proof of Lemma~1. 
\vskip 3pt
If a delay of the first type of duration $k_{i+1}^{(1)}$ for the cluster $i+1$ ends at time  $t$
$(1\le k_{i+1}^1\le l),$ then
$x_i=l-k_{i+1}^{(1)},$ $x_{i+1}=m,$
$\Delta(t)=m-l+k_{i+1}^{(1)},$ $i=0,1,\dots,N-1.$   

If, at time $t,$ a delay of the second type for the cluster $i$ of duration $k_{i}^{(2)}$ $(1\le k_i^{(2)}\le l-1)$ ends, then 
$x_i=0,$ $x_{i+1}=m+l-k_i^{(2)},$ 
$\Delta(t)=m+l-k_{i+1}^{(2)},$ $i=0,1,\dots,N-1.$    

Suppose, at time $t_0,$ a delay of the first type class for the cluster $i+1$ of duration
$k_i^{(1)}.$ When $x_i(t_0)=l-k_{i+1}^{(1)},$ $x_{i+1}(t_0)=m,$ $\Delta(t_0)=m-l+k_i^{(1)},$ $1\le k_{i+1}^{(1)}\le l),$ $i=0,1,\dots,N-1.$ 

Let us consider the joint behavior of the contours $i$ and $i+1$ ti the moment $t_0.$ 

For the moment $t_1=t_0-a_1,$ where $a_1=m+l-k_{i+1},$ the condition $x_i(t_1)=m$ holds, and, if, at time $t_1,$ a delay of the cluster $i,$ then $t_1$ is the moment of this delay end. Denote by $k_i^{(1)},$ the duration of this delay. If there was no delay,  then we assume that 
$k_i^{(1)}=0.$ We have $0\le k_i^{(1)}\le l.$ The delay that ends at time  $t_1,$ begins at time  
$t_2=t_1-a_1-k_i^{(1)},$ where
$$a_1+k_i^{(1)}=m+l+(k_i^{(1)}-k_{i+1}^{(1)}).\eqno(2)$$

If, at time $t_2=t_0-m,$ at the node $(i+1,i+2)$ such that $x_{i-1}=0,$ a delay of the cluster $C_{i+1}$ of duration $k_{i+1}^{(2)}.$ If there was no delay, then we assume  $k_{i+1}^{(2)}=0.$ We have $0\le k_{i+1}^{(2)}\le l-1.$ 
\vskip 3pt
{\bf Theorem 2.} {\it Under the conditions formulated above, then the values of delays  $k_i^{(1)}$ and $k_{i+1}^{(2)}$  satisfy the relation
$$k_i^{(1)}-k_{i+1}^{(2)}\ge k_{i+1}^{(1)}.\eqno(3)$$ 
\vskip 3pt
Proof.} Under the conditions of the theorem, the leading particle of the cluster $i$ is in the cell~1 (the contour $i$ is in the state 1) at time $t_0-a,$ where $a$ belongs to the set $I_1,$ which is the set $$[2m+k_i^1-k_{i+1}^{(1)},\ 2m+k_i^1-k_{i+1}^1+l-1],$$  
and the cluster $i+1$ occupies $m+1$ 
at time $t_0-a,$ where $a$ belongs to the set $I_2,$ which is the segment
$$[2m+k^{(2)}_{i+1}-l,\ 2m+k_{i+1}^{(2)}-1].$$  

Let us find the conditions for the values of delays such that the sets $I_1$ and $I_2$ do not intersect. 

Note first that the maximum element of the set $A_1$ is greater than $2m-1,$ and the minimum element of the set $A_2$ is not greater than $2m-1,$ and therefore the sets   $A_1$ and $A_2$ do not intersect if the minimum value of the set $I_1$ is greater than the maximum value of $I_2,$ i.e., if the condition 
$2m+k_i^{(1)}-k_{i+1}^{(1)}>2m+k_{i+1}^2,$ which is equivalent the inequality  (3) holds.        

Now we assume that the sets $I_1$ and $I_2$ intersect.  

Let $a_1$ be the maximum of the following elements: the minimum element of the set $I_1$ and the minimum element of the element $I_2,$ i.e.,
$$a_1=\max(2m+k_i^1-k^1_{i+1},2m+k_{i+1}^2-l).$$  
 At time $t_0-a_1,$ the leading particles of the clusters $i$ and $i+1$ are in the cells 1 and $m+1$ respectively. Hence, both the clusters have shifted onto at least one cell after their delays. Therefore, for  time $t_0-a_1-1,$ the delays of the clusters $i$ and $i+1$ have been ended, and the difference of the coordinates is equal to
$$\Delta_i(t_0-a_1)=m-l+k_{i+1}^{(1)})-k_i^{(1)}+k_{i+1}^{(2)}.$$  

A contradiction is that the system results in the non-admissible state. This contradiction is explained by that there is a delay of the cluster $i,$ $i+1,$ or the cluster that comes to the node later. 

Hence, in accordance with Lemmas 1 and 2  $\Delta_i(t_0-a_1)=m-l$ or $\Delta_i(t_0-a_1)=m+l$. 

In the first case, from $\Delta_i(t_0-a_1)=m-l,$ it follows the equality 
$$m-l+k_{i+1}^{(1)}-k_i^{(1)}+k^2_{i+1}=m-l.$$ 
From this, we get the equality
$$k_{i+1}^{(2)}=k_i^{(1)}-k_{i+1}^{(1)},$$
which is a special case of (3). The set $I_2$
$$I_2:\ [2m+k_i^{(1)}-l,m+k_i^{(1)}-k_{i+1}^{(1)}-1]$$
is located to the left of the set 
$$I_1:\ [2m+k_i^{(1)}-k_{i+1}^{(1)},2m+k_i^{(1)}-k_{i+1}^{(1)}-1].$$

In the second case, $\Delta(t_0-a_1)=m+l$ we get  
$$2l+k_i^{(1)}=k_{i+1}^{(1)}.$$
Since
$$k_{i+1}^{(1)}+k_{i+1}^{(2)}\le 2l-1,$$
we get a contradiction, which completes the proof of the theorem.
\vskip 3pt
{\bf Collorary 1.} {\it If, at time $t_0,$ a delay of the first class for the cluster $i+1$ begins, then the previous delay begins 
$t_0-m-l-k_{i+1}^{(1)}-k_i^{(1)}$ earlier.} 
\vskip 3pt
{\bf Theorem 3.} {\it Suppose, for a limit cycle, at time $t_0,$ a delay of the first type for the cluster $i$ begins, $i=0,1,\dots,N-1.$ Then, at time $t_0+m+l,$ a delay of the first type for the cluster $i+1$ begins (addition modulo $N).$       
\vskip 3pt
Proof.} In accordance with Theorem 1, we may construct the chain of delays 
$$k_{i+1}^{(1)}\le k_i^{(1)}\le k_{i-1}^{(1)}\le \dots $$ 
Since, for any $j,$ $k_j'\le L,$ we have that, from an index, all values are the same. Since the chain of delays is closed, we have that the durations of the delays $k_i^{(1)}$ are the same.   From the inequality (2), it follows that, for any $j,$ the duration of the delays of the second type satisfy the condition
$$k_{i+1}^{(2)}\le k_j^{(1)}-k_{j+1}^{(1)}=0$$
and, for this chain of delays, no delays occur. 

Taking into account (2) and the values $k_i^{(1)}$ are the same for any $j,$ we get that the duration of the time interval between the starts of the neighboring delays in the chain is equal to 
$$(m+L+(k_j^{(1)}-k_{j+1}^{(1)}))-0=m+l.$$
Theorem 3 has been proved.    
 \vskip 3pt   
Let any delay of the second type begins at time $t_0.$ Then we have $x_i(t_0)=0$ $\Delta x_i(t_0)=m+l-k_i^2-1.$  

At time $t_0-m,$ a delay of the first type for the cluster $i$ ends. The duration of this delay is equal to $k_i^{(1)},$ $0\le k_i^{(1)}\le l.$ The equality $(k_i^{(1)}=0$ means that there was no delay. 
\vskip 3pt
At time $t_0-m-l+k_i^{(2)},$ $x_{i+1}=t_2=m+l-k_i^{(2)},$ a delay of the second type for the cluster $i+1$ ends. Suppose that  
$k_i^{(2)}$ is the duration of this delay, $0\le  k_i^{(2)} \le l-1$ (if there was no delay, then $k_i^{(2)}=0).$
\vskip 3pt
{\bf Theorem 4.} {\it Under the conditions formulated above, the values of delays $k_i^{(1)}$ and $k_{i+1}^{(2)}$  satisfy the inequality
$$k_{i+1}^{(1)}-k_i^{(2)}\ge k_i^{(2}.$$
}
The proof of Theorem 4 is the same as the proof of Theorem 2.
\vskip 3pt
{\bf Theorem 5.} {\it Suppose, for a limit cycle, at moment $t_0,$ a delay of the first type for the cluster $i,$ $i=0,1,\dots,N-1,$ begins. Then, at time $t_0+m+l,$ a delay of the first type for the cluster $i+1$ begins (addition modulo $N).$}       
\vskip 3pt      
The proof of Theorem 4 is the same as th proof of Theorem 2.
\vskip 3pt
{\bf Theorem 6.} {\it In a cycle, either only delays of the first type occur, or only type of the second type occur.  
\vskip 3pt
Proof.} For the consistency with the previous theorems, let us introduce time countdown. Assume that, at $t=0,$ a delay of the first type of the cluster $i$ begins, and, at $t=t_2,$ a delay of the second type for the cluster $j$ begins. Let us pass along the chain of delays with step $m+l$ and find the moment $t_3\in (0,m+l)$ from that a delay of the cluster $k$ begins. Suppose that   $k<n.$ 
Moving consistently among the two chains of delays, we get that, at time $t=N(m+l),$ the delay of the second type for the cluster $k+N$ begins. 

Consider first the case in that either the indexes $i$ and $k$ are even or the indexes $n$ and $k$ are odds. Then there exist $n$ and $z$ such that $z=i-n=n-k.$ Therefore, at time $t=n(m+l),$ a delay of the first type for the cluster $i$ begins, at time $t=(n+1)(m+l),$ a delay of the first type for the cluster $z-1$ begins, and, at time $t=n(m+l)+t_3,$ a delay of the cluster $z$ begins. This  contradicts Theorem~2. 

Consider the case in that either the index $i$ is even and the index $k$ is odd or the index $i$ is odd and the index $k$ is even. In this case, the values $n$ and $z$ are prescribed such as $z=i+n-1=n-k.$ Then, at time $t=n(m+l)+t_3,$ a delay of the second type for the cluster $z$ begins, at time $t=(n+1)(m+l)$ a delay of the first type for the cluster $z+1$ begins, and, at time $t=(n+1)(m+l),$ a delay of the first type for the cluster $i$ begins. This contradicts Theorem~5. This contradiction completes the proof of Theorem~6.   
 \vskip 3pt
The proved theorems allow to describe an algorithm to get initial states and cycles, and to obtain the dependence of the spectrum from the parameters of the system. 

Let us formulate a sequence of statement giving numerical results,
\vskip 3pt
1. The maximum value of a delay of the first type is equal to $l.$ The maximum value of the delay of the second type is equal to $l-1.$ There are no other differences, i.e., suppose that, during a delay, the leading particle is in the cell $m.$  
\vskip 3pt
2. If the duration of a cluster is equal to $k,$ then we have: 
\vskip 3pt
a) the duration of a full turn is equal to  
$2m+k$ steps; 
\vskip 3pt
b) and the delay of the previous cluster is also equal to $k$ but this delay occurred $m+l$ earlier. 
\vskip 3pt
3. Denote by $T$ the period, i.e., the duration of cycle. 
\vskip 3pt
a) Let $r$ be the number of turns of a cluster during a cycle; $k_i$ is the duration of delay during the $i$th turn. Then we have 
$$T=(2m)r+k_1+\dots+k_r,$$   
and the values $k_i$ satisfy the condition $0\le k_i\le  l$ for delays of the first type (consequence of the statement 2a).   
\vskip 3pt that 
b) $T=(m+l)n$ (consequence of the statement 2b).
\vskip 3pt
4. The values $r$ and $k_i,$ $i=1,\dots,r,$ may be found from the following condition
$$T=(m+l)N=2mr+k_1+\dots+k_r, \eqno(4)$$
following from the statement 3.
\vskip 3pt
5. Obviously, the average velocity is equal to
$$v=1-\frac{k_1+\dots+k_r}{T}.\eqno(5)$$ 
\vskip 3pt
6. The condition (4) is not only necessary but is also sufficient. The initial states of cycle is formed follows.
First the sequence of states of a cluster. The begin of the cluster is in the cell $m-1$ (or $2m-1)$ $k_i+1$ times in the $i$th turn. The same sequence with the shift, in time, onto $m+l$ units gives the states of the next (or previous) clusters, and we must take into account the period $T.$ 
\vskip 3pt
{\bf Example 1.} Suppose $m=5,$ $n=3,$ $l=2.$

From (4) we get
$$(5+2)\times 3=10r+k_1+k_2.$$     
From this, we get $r=2,$ $k_1+k_2=1.$ Suppose $k_1=0,$ $k_2=1.$ We have a sequence of states of the contour 0
$$(0;1;2;3;4;5;6;7;8;9;0;1;2;3;4;4;5;6;7;8;9).$$
A shift onto $m+l$ positions to the left takes place. We have a sequence of states of the cluster 1
$$(4;4;5;6;7;8;9;0;1;2;3;4;5;6;7;8;9;0;1;2;3).$$
We have a sequence of states of the cluster 2
$$(7;8;9;0;1;2;3;4;4;5;6;7;8;9;0;1;2;3;4;5;6)$$
We return to the cluster~0:       
$$(0;1;2;3;4;5;6;7;8;9;0;1;2;3;4;4;5;6;7;8;9).$$
We have : $(0,4,7).$

There exists the following cycle, which is not got from the previous cycle by a permutation of the initial state) is
$(7,4,0).$ For this cycle, the delays are delays are delays of the second type. 

The average velocity of clusters equals
$$v=1-\frac{1}{21}=\frac{20}{21}.$$

There are two states depending on the initial state. They are the value 1 (the state of free movement, for example, the state $(0,0,0))$ and the value 20/21. 

\section{Set of values of average velocities of clusters}
 
 This section proves theorems regarding the values of the average velocities of clusters.
\vskip 7pt
{\bf Theorem 7.} {\it A necessary and conditional condition for a state of free movement to exist is 
$$l\le m.\eqno(6)$$ 
\vskip 3pt
Proof.} If the inequality (6) holds, then, for example, the state $(0,0,\dots,0)$ is a state of free movement.

If the inequality (6) does not hold, then the system results in the state of free movement from any initial state. Theorem~7 have been proved.
\vskip 3pt
{\bf Theorem 8.} {\it If the inequality (6) holds, then, for the limit cycle, the average velocity of clusters is not less than 2/3.  Moreover, if the inequality $l<m,$ then the average velocity greater 2/3.
\vskip 3pt
Proof.} Taking into account (4)--(6) that the duration of each delay of cluster is not greater than $l,$ we get Theorem~8. 
\vskip 3pt
{\bf Theorem 9.} {\it For prescribed values of $N,$ $m,$ the spectrum contains no more than $[N/3]$ values ($[a]$ is the integral part of the number $a).$
\vskip 3pt
Proof.} If the average velocity $v$ is less than 1, then, according to (4), (5) and Theorem~8,  for any initial state, the average velocity of clusters $v$ satisfies the condition
$$\frac{2}{3}\le v=\frac{2mr}{N(m+l)}<1,\eqno(7)$$
where $r$ is a natural number. From (7), we get
$$\frac{N}{3}+\frac{Nl}{3m}\le r=\frac{2m}{N(m+l)}<1.\eqno(8).$$
The number of integer values of $r$ satisfying (8) is not greater than   
$$\frac{2m}{N(m+l)}-\left(\frac{N}{3}+\frac{Nl}{3m}\right)+1=$$
$$=\frac{N}{6}+\frac{Nl}{3}+1.\eqno(9)$$
Combining (5) and (9), we get Theorem~8.
\vskip 3pt
{\bf Example 2.} Suppose $m=1$ a ({\it binary chain of contours}). Then according the results of [18], for any $k=0,1,\dots,[N/3],$ one may prescribe a state such that this state belongs to a limit cycle with average velocity   
$$v=1-k/N.\eqno(10)$$
For example, the state may be prescribed such that there are $k$ particles with cells with index~0, and, between any particles occupying  cells with index~0, there are no less than two particles occupying cells with index~1. There are not exist any limit cycles that do not satisfy (10). Depending on the initial state, the average velocity of cluster is equal to one of $[N/3]+1$ values, and $[N/3]$ of these values are less than~1.  
\vskip 3pt
{\bf Theorem 10.} {\it If $N$ is odd and the inequality 
$$l\le \frac{m}{N},(11)$$  
holds, then the system results in the state of free movement from any initial state.
\vskip 3pt
Proof.} For any initial state, it follows from (4) that 
$$r=\frac{N}{2}+\frac{Nl}{2m}-
\frac{k_1+\dots+k_r}{2m}.\eqno{12}$$ 
For at least one value $i,$ the value $k_i$ is positive. From this, taking into account the equality $0\le k_i\le l,$ $i=1,\dots,r,$ we get that, if (11) holds, the right hand of (12) can not be any integer number. Since the number $r$ is integer, the equality (12) does not hold, no limit cycles with the average velocity less than 1. From this, Theorem~10 holds.   
\vskip 3pt     
 {\bf Theorem 11.} {\it If the number $N$ is odd and 
$$l\le \frac{m}{N},(13)$$  
then the system results in a state of free movement from any initial state.}
\vskip 3pt
The proof of Theorem 11 as the proof of Theorem~10. 
\vskip 3pt 
{\bf Theorem 12.} {\it  
1) If the number $N$ is odd and the inequality
$$l\le \frac{m}{N}$$
holds, or the number $N$ is even and the inequality 
$$l\le \frac{2m}{N}$$
holds, then the system results in a state of free movement from any initial state. 
\vskip 3pt
2) If the number $N$ is odd and the condition 
$$\frac{m}{N}<l\le \frac{N}{2}$$ 
holds, or the number $N$ is even and the inequality
$$\frac{m}{N}<l\le \frac{N}{2},$$
holds, then limit cycles exist with the average equal to 1, and limit cycles may exist such that the average velocity is equal to one of $[N/3]$ values less  than~1. These values belong to the set of values computed according the formula 
$$v=\frac{2mr}{N(m+l)},$$  
where $r$ is an integer number satisfying the condition 
$2/3\le v<1.$ There exist a limit cycle satisfying this condition if and only if there exist values   $k_1,\dots,k_r$ satisfying (4).   
\vskip 3pt
3) If $l> m,$ then the system results in the state of collapse from any initial state. 
}
\vskip 3pt
Theorem 12 follows from (4), (5) and Theorems 1, 7, 9--11.

\section{Conclusion} 

This paper studies a discrete dynamical system. This system belongs to the class of the Buslaev contour networks. The studied system is called a closed chain of contours. 

The system contains $N$ closed sequences of cells. There are $2m$ cells and $l$ particles such that, at each discrete moment, these particles are located in adjacent cells and move simultaneously under a prescribed rule. 

We have proved that, if the number $N$ is odd and $l<m/N,$ or  
the number $N$ is even and $l<2m/N,$ then the system results in a state of free movement from any initial state. If the number $N$ is odd and $m/N\le l<N/2,$ or the number $N$ is even and $2m/N\le l<N/2,$ then there exist limit cycles with the average velocity equal to 1 and limit cycles may exist such that, depending on the initial state, the average velocity is equal to one of no more than  $[N/3]$ values less than~1. We have obtained approach to find these values of the average velocities and limit cycles corresponding to these velocities.

\end{document}